\documentclass{article}
\pdfoutput=1

\usepackage[title]{appendix}

\usepackage{amsthm, amssymb,amsmath}
\usepackage{url}\usepackage{hyperref}
\usepackage{authblk}

\newtheorem{theorem}{Theorem}
\newtheorem{corollary}{Corollary}
\newtheorem{proposition}{Proposition}
\newtheorem{remark}{Remark}
\newtheorem{lemma}{Lemma}
\theoremstyle{definition}
\newtheorem{example}{Example}
\theoremstyle{definition}
\newtheorem{definition}{Definition}

\begin{document}

	\title{Linear Differential Equation with Formal Power Series Non-Homogeneity Over a Ring with a Non-Archimedean Valuation
	\let\thefootnote\relax\footnote{The research was supported  by the  National Research Foundation of Ukraine funded by Ukrainian State budget in frames of project 2020.02/0096 ``Operators in infinite-dimensional spaces:  the interplay between geometry, algebra and topology''.}}

\author{Sergey Gefter\thanks{gefter@karazin.ua} }
\author{Anna Goncharuk\thanks{angoncharuk@ukr.net} }
\affil{Department of Mathematics \& Computer Sciences\\
	 V. N. Karazin Kharkiv National University}

	\date{}
		\maketitle
		
		\begin{abstract}
			Consider the linear differential equation of $m$-th order with constant coefficients from the valuation ring $K$ of a non-Archimedean field. We get sufficient conditions of uniqueness and existence for the solution of this equation from $K[[x]]$. Also the fundamental solution from $\frac{1}{x}K[[\frac{1}{x}]]$ of the equation is obtained and it is shown that the convolution of the fundamental solution and a non-homogeneity is a unique solution of the equation.
			
		\end{abstract}

\section{Introduction}
Let us consider the linear differential equation with constant coefficients
\begin{equation}\label{EQ}
	a_mw^{(m)}(x)+a_{m-1}w^{(m-1)}(x)+\ldots+a_1w'(x)+a_0w(x)=f(x),\  a_m\neq 0.
\end{equation}

Suppose we consider this equation on $\mathbb{R}$ and $f(x)$ is a continuous function. Then the Cauchy problem for Equation \eqref{EQ} with the initial condition $w(0)=w_0^0, w'(0)=w_0^1, \ldots, w^{(m-1)}(0)=w_0^{m-1}$ has a unique solution
\begin{equation}\label{Cauchy_nh}
	w(x)=w_h(x)+\int_0^x K(x,\xi)f(\xi)\ d\xi,
\end{equation}
where $w_h(x)$ is a solution of the homogeneous equation with these initial conditions and $K(x,\xi)$ is the Cauchy function of Equation \eqref{EQ}.

Suppose now $f(x)$ is a polynomial and $a_0\ne 0$. Then it is well known, that \eqref{EQ} has a unique polynomial solution. Despite this, it is difficult to find it using \eqref{Cauchy_nh}, because it is not known what initial condition this solution corresponds.

On the other hand, the algebraic method of undetermined coefficients allow us to find the unique polynomial solutions of Equation \eqref{EQ} if  we consider it in the ring $F[x]$, where $F$ is a field of characteristic zero and $a_0\neq 0$.

If now coefficients of $f(x)$ belong to an integral domain $K$, then the method of undetermined coefficients gives us a solution with coefficients from the quotient field of $K$. These coefficients may not belong to $K$ (see Example \ref{ex_n_r}).

We are studying Equation \eqref{EQ} with formal power series non-homogeneity. So let $f(x)$ be a formal power series with coefficients from the field $F$. Then the method of undetermined coefficients is useful for the Cauchy problem (\cite[Ch. VII]{Cartan}). It gives us infinitely many solutions: one for each initial value.

If $f(x)$ has coefficients from an integral domain $K$, in many interesting cases Equation \eqref{EQ} has a unique formal power series solution over $K$ (see Theorem \ref{uniq-nonarch}) and it is not known what initial value it corresponds to. So the method of undetermined coefficients is not useful in this case.

In 1808 B. Brisson \cite{Brisson} showed, particularly, that if denote $D=\frac{d}{dx}$, than the equation $Dy+y=f$ formally has a solution $y=f+Df+D^2f+\ldots$ (see also \cite{Pincherlet}). Note that using this formula for a polynomial non-homogeneity over a ring, one get the polynomial solution over that ring.

This result was generalized by U. Broggi for the case of Equation \eqref{EQ} (see \cite[\S5, 22.1]{Kamke}). He searched for the solution of \eqref{EQ} as a series
\begin{equation}\label{Br}
	\sum_{k=0}^{\infty} c_kf^{(k)}(x),
\end{equation}
where coefficients $c_k$ satisfies the equality
$$	(a_ms^m+a_{m-1}s^{m-1}+\ldots+a_1s+a_0)^{-1}=c_0+c_1s+c_2s^2+c_3s^3+\ldots.
$$
All these constructions was considered in the classical situation, that is over the field of real or complex numbers only if some convergence conditions hold. For instance, if $f(x)=e^x$, than series \eqref{Br} does not converge. About convergence conditions of the series \eqref{Br} see, for example,  \cite[Chapter 3, Section 3]{Leont'ev} and \cite[Chapter 1]{Sikkema} .

In the present paper we consider another situation, where we can use this construction. In Section 2 a formal explanation of what one can mean saying ``the series \eqref{Br} formally satisfies Equation (1)'' is given. In Section 3 the case of polynomial non-homogeneity over an arbitrary integral domain $K$ is described completely. If non-homogeneity $f(x)$ is a formal power series with the coefficients from $K$ and it is not a polynomial, than the series \eqref{Br} is not well-defined since infinite sums of the ring elements are appeared as coefficients of that series (Theorem \ref{f-conv-pol}).

In \cite[Example 2.1]{Fund-sol} it is shown, that the equation $y'+y=1+x+x^2+x^3+\ldots$ has no solution from $\mathbb{Z}[[x]]$. Theorem 2.1 in \cite{Hurwitz} yields, that one can make sense to the coefficients of the formal solution \eqref{Br} of this equation by considering infinite sums of integer in the ring of $p$-adic integers $\mathbb{Z}_p$. Note, that $\mathbb{Z}_p$ is a valuation ring of the field $p$-adic numbers $\mathbb{Q}_p$ with respect to the standard non-Archimedean valuation (see, for example, \cite[Section 1.2]{Perez-Garcia}).

Suppose now $(F,|\cdot|)$ is an arbitrary field $F$ of characteristics zero with a non-Archemedean valuation $|\cdot|$ (\cite[Section 1.2]{Perez-Garcia}), and $K$ is its valuation ring: $K=\{s\in F:|s|\leq 1\}$ (\cite[Ch.XII, \S 4]{Lang}).

In this paper the sufficient conditions for the uniqueness of the solution of Equation \eqref{EQ} are given (Theorem \ref{uniq-nonarch}). If additionally  $F$ is complete with respect to $|\cdot|$, then the solution of Equation \eqref{EQ} exists (Theorem \ref{exist-nonarch}). Moreover, the solution has the form \eqref{Br}, which converges with respect to the topology of coefficientwise convergence (about this topology see \cite[Chapter 1, Section 3]{Grauert}).
This result is specified to the case $K=\mathbb{Z}_p$ (Theorem \ref{Z}).

In the article \cite{Hurwitz} we construct an analogue of the convolution (Hurwitz product) of the Euler series $\mathcal{E}_b(x)=\frac{1}{x}-\frac{1!b}{x^2}+\frac{2!b^2}{x^3}-\ldots$ and an arbitrary formal power series with integer coefficients. It is shown that the Euler series can be regard as the fundamental solution to the first order equation $by'(x)+f(x)=y(x)$. In Section 5 of the present paper we generalize this construction for the case of Equation \eqref{EQ} and arbitrary valuation ring of a complete field with respect to the non-Archimedean valuation. Also we write down the fundamental solution explicitly for the second order equation (Example \ref{ex_2_E}).

Finally note, that in the present paper only the linear differential equations in rings of formal power series are considered. About differential equations in spaces of converging $p$-adic power series, and equations over fields of positive characteristics see, for example \cite{Dwork2}, \cite{Dwork}, \cite{Khrennikov}, \cite{Kochubei} and \cite{Robba}.

\section{Formal construction for solutions}		

Let $K$ be an arbitrary unitary commutative ring. Consider  the ring $ K[[x]][[y]]$ of power series having the form
$$w(x,y)=\sum_{k=0}^{\infty} {w_k(x)y^k},$$
where $w_k\in K[[x]].$
Let the $K$-linear operator $\tilde{D}: K[[x]][[y]]\rightarrow  K[[x]][[y]]$ take each $w(x,y)$ to $w_x'(x,y)\cdot y$.
Let $f(x)\in K[[x]]$, then denote $\tilde{f}(x,y)=f(x)$. For $a_0,a_1,\ldots,a_m\in K$ consider the equation in $K[[x]][[y]]$
\begin{equation}\label{EQ-abs}
	a_m\tilde{D}^{m}w+a_{m-1}\tilde{D}^{m-1}w+\ldots+a_1\tilde{D}w+a_0w=\tilde{f}.
\end{equation}

\begin{theorem}[existence and uniqueness of a solution]\label{uniq-ex-abs}
Let $a_0$ be invertible. Then there exists a unique solution of Equation \eqref{EQ-abs} from $K[[x]][[y]]$:
\begin{equation}\label{sol_abs}
w(x,y)=\sum_{k=0}^{\infty} {c_kf^{(k)}(x)y^k},
\end{equation}
where coefficients $c_k$ are found from the equality
\begin{equation}\label{coeffs}
	(a_mt^m+a_{m-1}t^{m-1}+\ldots+a_1t+a_0)^{-1}=c_0+c_1t+c_2t^2+c_3t^3+\ldots.
\end{equation}
\end{theorem}

\begin{proof}
Since $a_0$ is invertible, the polynomial $P(t)=a_mt^m+a_{m-1}t^{m-1}+\ldots+a_1t+a_0$ is invertible in the ring $K[[x]]$. Equation \eqref{EQ-abs} can be written as $P(\tilde{D})w=\tilde{f}$. It has a unique solution of the form $(P(\tilde{D}))^{-1}\tilde{f}$, where $(P(t))^{-1}=\sum_{k=0}^{\infty} {c_kt^k}$. Its coefficients are uniquely determined by \eqref{coeffs}. Since $(\tilde{D}^k\tilde{f})(x,y)=f^{(k)}(x)y^k$, we obtain the solution \eqref{sol_abs}.

\end{proof}

\begin{remark}\label{low-order-c}
Let us obtain an explicit formula for $c_i$ for the first and second order equations.

In the case $n=1$ we have
$$\frac{1}{a_0+a_1t}=\sum_{j=0}^{\infty} {(-1)^ja_0^{-j-1}a_1^jt^j}.$$
Thus, $c_j=(-1)^ja_0^{-j-1}a_1^j$.

In the case $n=2$ we have
\begin{multline*}
	\frac{1}{a_0+a_1t+a_2t^2}=\sum_{n=0}^{\infty} {(-1)^na_0^{-n-1}\sum_{k=0}^{n} {\binom{n}{k}a_1^{n-k}a_2^k}t^{n+k}}=\\
	=\sum_{k=0}^{\infty} \sum_{n=k}^{\infty} {(-1)^na_0^{-n-1} {\binom{n}{k}a_1^{n-k}a_2^k}t^{n+k}}=\\
	=\sum_{k=0}^{\infty} \sum_{n=0}^{\infty} {(-1)^{n+k}a_0^{-n-k-1} {\frac{(n+k)!}{k!n!}a_1^{n}a_2^k}t^{n+2k}}=\\
	=\sum_{j=0}^{\infty} \left( \sum_{k=0}^{[\frac{j}{2}]} { {\frac{(-1)^{j-k}(j-k)!}{k!(j-2k)!}a_0^{-j+k-1}a_1^{j-2k}a_2^k}}\right)t^{j}.
\end{multline*}
Therefore,
\begin{multline*}
c_j=\sum_{k=0}^{[\frac{j}{2}]} { {\frac{(-1)^{j-k}(j-k)!}{k!(j-2k)!}a_0^{-j+k-1}a_1^{j-2k}a_2^k}}=\\
=\sum_{k=0}^{[\frac{j}{2}]} { {(-1)^{j-k}\binom{j-k}{k}a_0^{-j+k-1}a_1^{j-2k}a_2^k}}.
\end{multline*}
\end{remark}

Let us consider the differential equation \eqref{EQ} on $K$ with the non-homogeneity $f(x) \in K[[x]]$ and find a connection between Equation \eqref{EQ} and Equation \eqref{EQ-abs}.

Obviously, if $f(x)$ is a polynomial, then \eqref{Br} solves \eqref{EQ}. On the uniqueness of such solution see Proposition \ref{uniq-p} in Section 3.

\begin{definition}
	Let $w_k(x) \in K[[x]]$. We say that the series $w(x)=\sum_{k=0}^{\infty} w_k(x)$ formally satisfies Equation \eqref{EQ} if the series $w(x,y)=\sum_{k=0}^{\infty} w_k(x)y^k$ satisfies Equation \eqref{EQ-abs}.
\end{definition}

\begin{theorem}\label{formal}
	Suppose the coefficient $a_0$ of Equation \eqref{EQ} is invertible.
	\begin{itemize}
		\item If $\{c_k\}$ satisfies \eqref{coeffs}, then the series $\sum_{k=0}^{\infty} {c_kf^{(k)}(x)}$ formally satisfies \eqref{EQ}.
		\item Let $K$ be an integral domain and torsion-free $\mathbb{Z}$-module, i.e. for $n\in \mathbb{Z}\setminus\{0\}$ and $a\in R$ the equality $na=0$ implies $a=0$ (see \cite[Chapter VII, \S 2]{Bourbaki}) and let $f\in K[[x]]\setminus K[x]$. If $\sum_{k=0}^{\infty} {c_kf^{(k)}(x)}$ formally satisfies \eqref{EQ}, then $\{c_k\}$ satisfies the equality \eqref{coeffs}.
	\end{itemize}
\begin{proof}
Theorem \ref{uniq-ex-abs} directly yields the first statement.

Let now $b_k\in K,\ k=1,2,3,\ldots$ and the series  $\sum_{k=0}^{\infty} {b_kf^{(k)}(x)}$ satisfies Equation \eqref{EQ}. If $c_k$ are found from \eqref{coeffs}, then the series $\sum_{k=0}^{\infty} {(b_k-c_k)f^{(k)}(x)}y^k$ is a solution of the homogeneous equation 
\begin{equation}
	a_m\tilde{D}^{m}w+a_{m-1}\tilde{D}^{m-1}w+\ldots+a_1\tilde{D}w+a_0w=0.
\end{equation}
By Theorem \ref{uniq-ex-abs}, it follows that $(b_k-c_k)f^{(k)}(x)=0$ for any $k$.

Let $b_j\neq c_j$ for some $j$. Since $K$ is integral domain, then $\frac{(j+i)!}{j!}f_{j+i}=0$ for any $i$. Since $K$ is a torsion-free $\mathbb{Z}$-module, we obtain that $f_j=f_{j+1}=\ldots=0$, i.e. $f\in K[x]$. The contradiction concludes the proof.
\end{proof}

\end{theorem}

Let us show that if the series $\sum_{k=0}^{\infty} w_k(x)$, where $w_k\in K[[x]]$, formally satisfies Equation \eqref{EQ} and is well-defined, then the sum of this series is a solution of Equation \eqref{EQ}. To this end we need the following lemma.

\begin{lemma}\label{S}
Suppose some $K$-linear map $S:\ K[[x]][[y]] \rightarrow K[[x]]$ with a domain $\mathfrak{D}(S)$ satisfies the following conditions:
\begin{itemize}
	\item $\tilde{D}(\mathfrak{D}(S))\subset \mathfrak{D}(S)$;
	\item $S(\tilde{D}w) =(Sw)'(x)$ for any $w\in \mathfrak{D}$;
	\item for any $f\in K[[x]]$, if $\tilde{f}(x,y)=f(x)$ then $S\tilde{f}=f$.
\end{itemize}
Let $w \in \mathfrak{D}(S)$ and $w$ solves \eqref{EQ-abs}. Then $Sw$ solves \eqref{EQ}.
\end{lemma}

\begin{proof}
Since $w$ belongs to $\mathfrak{D}(S)$, then $\tilde{D}w, \tilde{D}^2w, \ldots, \tilde{D}^mw$ belong to $\mathfrak{D}(S)$.
	
	The series $w(x,y)\in \mathfrak{D}(S)$ is a solution of \eqref{EQ-abs}, so $$S(a_m\tilde{D}^{m}w+a_{m-1}\tilde{D}^{m-1}w+\ldots+a_1\tilde{D}w+a_0w)=S\tilde{f}.$$
	
	Since for any $j$ we have
	$S(\tilde{D}^{j}w(x,y))=(Sw)^{(j)}$ and $S\tilde{f}=f$, then
	$$a_m(Sw)^{(m)}+a_{m-1}(Sw)^{(m-1)}+\ldots+a_1(Sw)'+a_0Sw=f.$$
\end{proof}

In the following proposition we consider the map $S: K[[x]][[y]] \rightarrow K[[x]]$ be given by
\begin{equation}\label{map_s}
v(x,y) \mapsto v(x,1).
\end{equation}

\begin{corollary}\label{connect}
Suppose the series $w(x)=\sum_{k=0}^{\infty} {w_k(x)}$ converges in the ring $K[[x]]$ in the Krull topology \cite[Chapter 1, Section 3]{Grauert}. If it formally satisfies \eqref{EQ}, then it is a solution of \eqref{EQ}.
In this case the series $\sum_{n=0}^{\infty} c_n \binom{n+k}{n} n!f_{k+n}$ has a finite number of terms for any $k$ and
\begin{equation}\label{sol}
	w(x)=\sum_{k=0}^{\infty} {\left(\sum_{n=0}^{\infty} c_n \binom{n+k}{n} n!f_{k+n}\right)x^k }.
\end{equation}
It can also be written as
\begin{equation}\label{sol-f}
	w(x)=\sum_{k=0}^{\infty} {c_kf^{(k)}(x)}.
\end{equation}
\end{corollary}

\begin{proof}
	Let the map $S:\ K[[x]][[y]] \rightarrow K[[x]]$ with the domain  
	$$\mathfrak{D}(S)=\{v(x,y)=\sum_{k=0}^{\infty} v_k(x)y^k\ :\ \sum_{k=0}^{\infty} {v_k(x)} \mbox{ converges in the Krull topology}\}$$
	be given by \eqref{map_s}.
Note that the series $\sum_{k=0}^{\infty} {v_k(x)}$, where $v_k(x)=\sum_{j=0} v_{kj}x^j$ converges in the Krull topology if for any $j$ series $\sum_{k=0}^{\infty} v_{kj}$ have a finite number of terms.
	
The map $S$ satisfies the conditions of Lemma \ref{S}. Indeed, $\sum_{k=0}^{\infty} {v_k(x)}$ converges in the Krull topology, then also $\sum_{k=0}^{\infty} {v'_k(x)}$  converges in the Krull topology, thus $v'_x(x,y)y\in \mathfrak{D}(S)$. Also $S\tilde{f}(x,y)=f(x)$ for $f\in K[[x]]$ and $S(\tilde{D}w)=S(w'_x(x,y)y)=(w'_x(x,y)y)\Big|_{y=1}=w'_x(x,1)=(Sw)'$.
	
Let $w(x,y)=\sum_{k=0}^{\infty} {w_k(x)y^k}$. Then by Lemma \ref{S} $w(x)=\sum_{k=0}^{\infty} {w_k(x)}$	solves \eqref{EQ}. By Theorem \ref{uniq-ex-abs}
 $w(x)$ can be written as \eqref{sol} and \eqref{sol-f}.
\end{proof}

\begin{remark}
	Further it will be shown for the integral domain $K$ that the series \eqref{sol-f} converges in the Krull topology if and only if $f$ is a polynomial (see Theorem \ref{f-conv-pol}).
\end{remark}

\section{Polynomial solution}

Suppose $K$ is an integral domain and $F=Frak(K)$ is the quotient field of $K$.

\begin{proposition}[existence and uniqueness of a polynomial solution]\label{uniq-p}
	Suppose $a_0 \neq 0$. Then the following statements hold:
\begin{itemize}
	\item Let $f(x)\in F[x]$. Then Equation \eqref{EQ} has a unique solution from $F[x]$. This solution has a form \eqref{sol-f}. 
	\item If $f(x)\in F[[x]]\setminus F[x]$, then Equation \eqref{EQ} has no solution from $F[x]$.
\end{itemize}
\end{proposition}

\begin{proof}
First let us prove that the Equation \eqref{EQ} has no more then one solution from $F[x]$. It is enough to prove that the homogeneous equation
\begin{equation}\label{EQ-hom}
	a_mw^{(m)}(x)+a_{m-1}w^{(m-1)}(x)+\ldots+a_1w'(x)+a_0w(x)=0,
\end{equation}	
has only zero polynomial solution.
Suppose the polynomial $w(x)=\sum_{i=0}^{n} w_ix^i$ is a non-zero solution of \eqref{EQ-hom} and $\deg{w}=n$. Then the coefficient of $x^n$ in the left-hand side of this equation is $a_0w_n=0$. Since $a_0\neq 0$, then $w_n=0$, which contradicts $\deg{w}=n$. Thus, \eqref{EQ} has no more than one solution from $F[x]$.

Since $f(x)$ is a polynomial, then left-hand side of \eqref{sol-f} is a polynomial too. By Corollary \ref{connect} it is a solution of Equation \eqref{EQ}.

Finally, if the polynomial $w(x)$ solves Equation \eqref{EQ}, then the non-homogeneity of \eqref{EQ} equals  $a_mw^{(m)}(x)+a_{m-1}w^{(m-1)}(x)+\ldots+a_1w'(x)+a_0w(x)$, which is a polynomial.
\end{proof}

If $f\in F[x]$ and $a_0\neq 0$, then the formula \eqref{sol-f} gives a solution from $F[x]$. Then the previous theorem yields the following result for the solution from $K[x]$.

\begin{corollary} Suppose the non-homogeneity of \eqref{EQ} belongs to $K[x]$.
	\begin{enumerate}
		\item If $a_0$ is invertible, then the solution \eqref{sol-f} of Equation \eqref{EQ} belongs to $K[x]$.
		\item If $a_0$ is not invertible, then there exists $f\in K[x]$ for which Equation \eqref{EQ} has no polynomial solution.
		\item If $a_0\neq 0$ is not invertible, then either the solution \eqref{sol-f} belongs to $K[x]$ or Equation \eqref{EQ} has no solution from $K[x]$.
		\item If $a_0=0$ and $K$ is infinite, then either there are infinitely many solutions of \eqref{EQ} from $K[x]$ or no one.
	\end{enumerate}
\end{corollary}

\begin{proof}
	Statements 1 and 3 directly follows by Proposition \ref{uniq-p}.
	
	To prove the statement 2 it is enough to consider $f(x)=1$. Then
	$$a_mw^{(m)}(x)+a_{m-1}w^{(m-1)}(x)+\ldots+a_1w'(x)+a_0w(x)=1.$$
	Then the maximal degree of $w$ equals to $0$. Therefore, the only solution of this equation is constant $w(x)=C$. We get $a_0C=1$. It means that $a_0$ is invertible.
		
	Now let us prove 4. Denote by $l$ the minimal number for that $a_l\neq 0$. Now set $v(x)=w^{(l)}(x)$. Then by Proposition \ref{uniq-p} there exists a unique $v_0\in F[x]$ that satisfies \eqref{EQ}. We get the equation $w^{(l)}(x)=v_0(x)$. If it has a solution from $K[x]$, then we get infinitely many solutions by adding constants from $K$.
\end{proof}

\begin{example}
	Let consider the equation $3w'+2w=2x+5$ over $\mathbb{Z}$ and $\mathbb{Q}$. In this case $a_0=2\neq 0,\ a_1=3$, then $c_0=\frac{1}{a_0}=\frac{1}{2}$ and $c_1=-\frac{a_1}{a_0^2}=-\frac{3}{4}$. Therefore, by formula~\eqref{sol-f}, the unique solution from $\mathbb{Q}[x]$ has the form
	$$w(x)=c_0f(x)+c_1f'(x)=\frac{1}{2}(2x+5)-2\cdot\frac{3}{4}=x+1\mbox{ and } w(x) \in \mathbb{Z}[x].$$
\end{example}	

\begin{example}\label{ex_n_r}
	Let consider the equation $3w'+2w=x+5$. The left side of this equation coincides with the previous example, thus $c_0=\frac{1}{2},\ c_1=-\frac{3}{4}$. By formula~\eqref{sol-f}, the solution has the form
	$$w(x)=c_0f(x)+c_1f'(x)=\frac{1}{2}(x+5)-\frac{3}{4}=\frac{1}{2}x+\frac{7}{4}.$$
	By Proposition \ref{uniq-p}, this solution is unique in $\mathbb{Q}[x]$. Since it does not belong to $\mathbb{Z}[x]$, then the equation $3w'+2w=x+5$ has no solution in $\mathbb{Z}[x]$.
\end{example}

\begin{example}
	The equation $w'(x)=x$ has no solution from $\mathbb{Z}[x]$ and it has infinitely many solutions from $\mathbb{Z}_{(3)}[x]$, where
	$$\mathbb{Z}_{(3)}=\{\frac{n}{m}:\ (n,m)=1 \mbox{ and } m \mbox{ is not divisible by } 3\}.$$
	Indeed, $w(x)=\frac{x^2}{2}+C$ and only they are solutions of this equation in $\mathbb{Q}[x]$, all of them belong to $\mathbb{Z}_{(3)}[x]$ and no one among them belongs to $\mathbb{Z}[x]$.
\end{example}

Since for the first and second order equation the explicit form of $c_j$ are found in Remark \ref{low-order-c}, then one can right down the explicit formula of the polynomial solution.
\begin{corollary} Let $a_0$ be invertible and $f\in K[x]$. Then
	\begin{enumerate}
		\item The equation $a_1w'(x)+a_0w(x)=f(x)$ has a unique polynomial solution
		\begin{multline*}
			w(x)=\sum_{k=0}^{\infty} {(-1)^ka_0^{-k-1}a_1^kf^{(k)}(x)}=\\
			=\sum_{k=0}^{\infty} {\left( \sum_{n=0}^{\infty} (-1)^na_0^{-n-1}a_1^n \binom{n+k}{n} n!f_{k+n}\right) x^k}.
		\end{multline*}
		\item The equation $a_2w''(x)+a_1w'(x)+a_0w(x)=f(x)$ has a unique polynomial solution
		\begin{multline*}
			w(x)=\sum_{j=0}^{\infty} {f^{(j)}(x)\sum_{k=0}^{[\frac{j}{2}]} { {(-1)^{j-k}\binom{j-k}{k}a_0^{-j+k-1}a_1^{j-2k}a_2^k}}}=\\
		=\sum_{k=0}^{\infty} \left({\sum_{j=0}^{\infty} \binom{j+k}{j} j!f_{k+j}}\sum_{n=0}^{[\frac{j}{2}]} { {(-1)^{j-n}\binom{j-n}{n}a_0^{-j+n-1}a_1^{j-2n}a_2^n}}\right) x^k.
	\end{multline*}
\end{enumerate}

\end{corollary}

Suppose now the characteristic of $F$ equals zero. Consider Equation \eqref{EQ} with a non-homogeneity from $F[[x]]$. By Theorem \ref{formal}, the series \eqref{sol-f} formally satisfies \eqref{EQ}. Theorem \ref{f-conv-pol} below shows that if $f$ is not a polynomial, then \eqref{sol-f} is not well-defined. First we prove the following condition of convergence.
 
\begin{lemma}
Let $c_j,f_j\in F,\ j=0,1,2,\ldots$ be arbitrary elements in $F$. Then the series \eqref{sol} converges in the Krull topology on $F[[x]]$ if and only if for any $k$ there exists $i$ such that $c_jf_{j+k}=0$ for any $j>i$.	
\end{lemma}

\begin{proof}
	Indeed, \eqref{sol} converges if and only if $\sum_{j=0}^{\infty} {c_j(k+j)!f_{k+j}}$ has a finite number of non-zero terms for any $k$.
\end{proof}

The following example shows that the condition ``for any $k$ there exists $i$ such that $c_jf_{j+k}=0$ for any $j>i$'' could hold even in the case both $\{c_j\}$ and $\{f_j\}$ are not finite. 

\begin{example}
	Let $F=\mathbb{Q}$, $c_j=1$ if there exist $r$ such that $j=2^r$. Otherwise set $c_i=0$. Set $f_i=1$ if  $i=2^r+r$ and $f_i=0$ otherwise.
	
	For any $k$ set $i=2^{k+1}$. Since for any $j>i$ such that $c_j\neq 0$ there exists $r$ such that $j=2^r$, then we have $j=2^r>i=2^{k+1}$. Therefore $k\leq r-1$. Then $f_{k+j}=f_{2^r+k}=0$, since $2^{r-1}+r-1<2^r<2^r+k<2^r+r$. We get that if $c_j\neq 0$, then $f_{k+j}=0$, so $c_jf_{j+k}=0$.
\end{example}

Despite this since $\{c_j\}$ is not an arbitrary sequence but constructed to satisfy the equality \eqref{coeffs}, then the following theorem holds.

\begin{theorem}\label{f-conv-pol}
	Suppose $f(x)\in F[[x]]$, $a_j \in F$, $a_0\neq 0$ and $c_j$ are found from the equality \eqref{coeffs}. Then the series \eqref{sol-f} converges in the Krull topology on $F[[x]]$ if and only if $f(x)\in F[x]$.
\end{theorem}

\begin{proof}
	Note that since $\{c_j\}$ satisfies the equality \eqref{coeffs}, it is not finite. Also it solves the system
	$$\left\{
		\begin{array}{rcl}
			a_0c_0&=&1,\\
			\sum_{i=0}^{j} {a_ic_{j-i}}&=&0,\ j=1,2,3,\ldots,m \\
			\sum_{i=0}^{m} {a_ic_{j-i}}&=&0,\ j=m+1,m+2,\ldots
		\end{array}
		\right.
	$$
	First note that for any $j$ there exists $1\leq i\leq m$ such that $c_{j+i}\neq 0$. Indeed, if there exists $j$ such that $c_{j+1}=c_{j+2}=\ldots=c_{j+m}=0$, then the system implies $c_i=0$ for any $i\geq j+1$. This contradicts the fact that $\{c_j\}$ is not finite.
		
	We claim that if $f(x)$ is not a polynomial, then the sum $\sum_{j=0}^{\infty} {c_j(k+j)!f_{k+j}}$ has an  infinite number of terms for some $k$ and then the series \eqref{sol} is not well-defined. Assume the converse. Then for any $k$ there exists $i_k$ such that for any $j>i_k$ the equality $c_jf_{k+j}=0$ holds. Consider $j> \max\limits_{k=0,1,\ldots,m} {i_k}$, such that $c_j\neq 0$. Then $f_{j}=f_{j+1}=\ldots=f_{j+m}=0$. As shown above, there exists $1\leq i\leq m$ such that $c_{j+i}\neq 0$, then $f_{j+i}=f_{j+i+1}=\ldots=f_{j+i+m}=0$. Therefore, since $j+i+m>j+m$, then for any $k>j$, $f_k=0$, that is $f(x)$ is a polynomial. This contradiction proves the theorem.
\end{proof}

\section{Ring with a non-Archimedean valuation}

\subsection{General non-Archimedean case}
Suppose $(F,|\cdot|)$ is a field with a non-Archimedean valuation and $K=\{s \in F: |s| \leq 1\}$ is its valuation ring. Note that an element $a\in K$ is invertible if and only if $|a|=1$. Now we consider the topology of coefficientwise convergence on $K[[x]]$  (see \cite[Chapter 1, Section 3]{Grauert}).

\begin{theorem}\label{uniq-nonarch}
 	Let in Equation \eqref{EQ} the coefficients $a_0,a_1,\ldots,a_m$ belong to $K$, $|a_0|=1$, $|a_i|<1$ for any $1\leq i\leq m$. Then this equation has no more than one solution in $K[[x]]$.
\end{theorem}
 
 \begin{proof}
 	Let us consider the homogeneous equation
 	\begin{equation}\label{hom}
 		a_my^{(m)}(x)+a_{m-1}y^{(m-1)}(x)+\ldots+a_1y'(x)+a_0y(x)=0.
 	\end{equation}
 	Find a solution of the form $y(x)=y_0+y_1x+y_2x^2+\ldots$.
 	Then for any $k$ the following equality holds:
 	
 	$$\sum_{i=0}^{m} {\frac{(k+i)!}{k!}a_iy_{k+i}}=0.$$
 	Then 
 	\begin{equation}\label{y_k}
 		y_k=-a_0^{-1}\sum_{i=0}^{m-1}\frac{(k+i+1)!}{k!}a_{i+1}y_{k+i+1}.
 	\end{equation}
  Now for any indexes $i_1,\ldots,i_k$ set  $s_k=\sum_{j=0}^{k} i_j$ and $p_k=\prod_{j=1}^{k} a_{i_j+1}$. We get
 	\begin{multline*}
 		y_0=-a_0^{-1}\sum_{i_1=0}^{m-1}(s_1+1)!p_1y_{s_1+1}=a_0^{-2}\sum_{i_1=0}^{m-1}p_1\sum_{i_2=0}^{m-1}(s_2+2)!a_{i_2+1}y_{s_2+2}=\\
 		=a_0^{-2}\sum_{i_1,i_2=0}^{m-1}(s_2+2)!p_2y_{s_2+2}=
 		-a_0^{-3}\sum_{i_1,i_2,i_3=0}^{m-1}(s_3+3)!p_3y_{s_3+3}=\ldots\\
 	\ldots=(-1)^{k}a_0^{-k}\sum_{i_1,i_2,\ldots,i_k=0}^{m-1}(s_k+k)!p_ky_{s_k+k}.
 	\end{multline*}
 	Set $b=\max_{i=1}^{m}{|a_i|}<1$. Then $|p_k|\leq b^k \rightarrow 0$ as $k\rightarrow \infty$. Since $|\cdot|$ is non-Archimedean, then, $|y_0|\leq b^k$, so $y_0=0$. Similarly, using \eqref{y_k} we get $y_k=0$ for any $k$.
 \end{proof}

Now let us obtain an existence condition for the solution. To this end, consider the map $S$ be given by \eqref{map_s}, which has the domain
\begin{multline*}
	\mathfrak{D}(S)=\{\sum_{i=0}^{\infty} {v_{i}(x)y^k}:\ \sum_{i=0}^{\infty} {v_{i}(x)} \mbox{ converges in the topology of } \\ 
	\mbox{coefficientwise convergence on } K[[x]]\}.
\end{multline*}

\begin{theorem}\label{exist-nonarch}
	Suppose $F$ is complete with respect to $|\cdot|$. Let $|a_0|=1$, $|a_i|<1$ for any $1\leq i\leq m$. Then the series \eqref{sol-f} converges by the topology of coefficientwise convergence in $K[[x]]$ and the sum of this series is a unique solution of \eqref{EQ} in $K[[x]]$.
\end{theorem}

\begin{proof}
We claim, that $S$ is under the conditions of Lemma \ref{S}. Indeed, similarly to Corollary \ref{connect}, it can be proven that $S(\tilde{D}v) =(Sv)'(x)$
and $S\tilde{f}(x,y)=f(x)$ for $f\in K[[x]]$. If $v\in \mathfrak{D}(S)$, then for any $j$ series $\sum_{i=0}^{\infty} {v_{ij}}$ converges in $K$ with respect to ${|\cdot|}$.
Then $\sum_{j=0}^{\infty} {(j+1)v_{i(j+1)}}$ converges in $K$ with respect to $|\cdot|$ for any $j$. Therefore $v'_x(x,y)=\sum_{i=0}^{\infty} \left( \sum_{j=0}^{\infty} {(j+1)v_{i(j+1)}x^{j}}\right) y^i \in \mathfrak{D}(S)$.

Let the series $w(x,y)$ be the solution \eqref{sol_abs} of Equation \eqref{EQ-abs}. By Lemma \ref{S}, if it  belongs to $\mathfrak{D}(S)$, then $Sw$ is a solution of \eqref{EQ}. Let us check that \eqref{sol_abs}
 belongs to $\mathfrak{D}(S)$. To this end, let us estimate the coefficients $c_j$ of $w(x,y)$. We get
$$\frac{1}{a_0+a_1t+a_2t^2+\ldots+a_mt^m}=\frac{a_0^{-1}}{1-a_0^{-1}(-a_1-a_2t-\ldots-a_mt^{m-1})t}=$$
$$=\sum_{n=0}^{\infty} {(-1)^na_0^{-n-1}(a_1+a_2t+\ldots+a_mt^{m-1})^nt^n}=$$
$$=\sum_{n=0}^{\infty} {\sum\limits_{i_1+i_2+\ldots+i_m=n} (-1)^n\frac{n!}{i_1!i_2!\cdot\ldots\cdot i_m!}\cdot \frac{a_1^{i_1}a_2^{i_2}\cdot\ldots\cdot a_m^{i_m}}{a_0^{n+1}} t^{n+i_2+2i_3+\ldots+(m-1)i_m}}$$

Since $n+i_2+2i_3+\ldots+(m-1)i_m\leq mn$, thus for the degree $j$ one can estimate $n \geq [\frac{j}{m}]$.

Since $|\cdot|$ is non-Archimedean, then $|c_j|$ is no more than maximal of all $$|(-1)^n\frac{n!}{i_1!i_2!\cdot\ldots\cdot i_m!}\cdot \frac{a_1^{i_1}a_2^{i_2}\cdot\ldots\cdot a_m^{i_m}}{a_0^{n+1}}|, \mbox{ such that } n+i_2+2i_3+\ldots+(m-1)i_m=j.$$
Set $b=\max_{i=1}^{m}{|a_i|}<1$. Then we get
\begin{equation}\label{ineq}
	|(-1)^na_0^{-n-1}\frac{n!}{i_1!i_2!\cdot\ldots\cdot i_m!}a_1^{i_1}a_2^{i_2}\cdot\ldots\cdot a_m^{i_m}|<b^n\leq b^{[\frac{j}{m}]}
\end{equation}
Since $b<1$, then
\begin{equation}\label{tends_c}
	c_j|\leq b^{[\frac{j}{m}]} \rightarrow 0 \mbox{ as } j\rightarrow\infty.
\end{equation}

Since $|\binom{j+k}{j} j!f_{k+j}|\leq 1$ for any $k$, then the series $\sum_{j=0}^{\infty} c_j \binom{j+k}{j} j!f_{k+j}$ from \eqref{sol} converges in $K$.  We used the fact that in a complete non-Archimedean field a series $\sum {a_k}$ converges if and only if the sequence $\{a_k\}$ tends to $0$ (\cite[Section 2.1]{Perez-Garcia}). Thus, $w\in \mathfrak{D}(S)$. The solution $Sw$ of \eqref{EQ} can be found as
$$w(x)=\sum_{k=0}^{\infty} {\left(\sum_{j=0}^{\infty} c_j \binom{j+k}{j} j!f_{k+j}\right)x^k }.$$
\end{proof}

\subsection{Case of integers}

In a particular case when $(F,|\cdot|)=(\mathbb{Q}_p,|\cdot|_p)$ is the field of $p$-adic numbers (\cite[Section 1.2]{Perez-Garcia}), then the valuation ring of $\mathbb{Q}_p$ is the ring of $p$-adic integers $\mathbb{Z}_p$.

The following result specifies Theorem \ref{exist-nonarch} in the case of the integer coefficients in Equation \eqref{EQ}.

\begin{theorem}\label{Z}
Suppose $a_0,a_1,\ldots,a_m$ are integer. Then Equation \eqref{EQ} has a unique solution from $\mathbb{Z}_p[[x]]$ for any prime $p$ that is not a divisor of $a_0$.
\end{theorem}

\begin{proof}
If $p$ does not divide $a_0$, then $|a_0|_p=1$. The series $\sum_{n=0}^{\infty} c_n \binom{n+k}{n} n!f_{k+n}$ converges in $\mathbb{Z}_p$, because $|c_n \binom{n+k}{n} f_{k+n}|_p\leq 1$ and $|n!|_p$ tends to $0$. Thus, using the same $S$ as in the previous theorem, we conclude that the solution from $\mathbb{Z}[[x]]$ of Equation \eqref{EQ} exists.

Now let us prove the uniqueness. For any indexes $i_1,\ldots,i_k$ set  $s_k=\sum_{j=0}^{k} i_j$. The coefficients of a solution $y(x)$ of the homogeneous equation satisfy the following equality
$$y_0=k!\cdot (-1)^{k}a_0^{-k}\sum_{i_k=0,k\in\mathbb{N}}^{m-1}\frac{(s_k+k)!}{k!}\prod_{j=1}^{k} a_{i_j+1}y_{s_k+k}$$
for any $k$. Since  
$|(-1)^{k}a_0^{-k}\sum_{i_k=0,k\in\mathbb{N}}^{m-1}\frac{(s_k+k)!}{k!}\prod_{j=1}^{k} a_{i_j+1}y_{s_k+k}|_p\leq 1$ and $|k!|_p\rightarrow 0$ as $k\rightarrow \infty$, then $y_0=0$. Similarly we get $y_k=0$ for any $k$.
\end{proof}

\begin{remark}
	This solution does not necessarily belong to $\mathbb{Z}[[x]]$. For example, the equation $y'+y=1+x+x^2+x^3+\ldots$ due to the previous theorem has a unique solution in $\mathbb{Z}_p[[x]]$ and, as it is stated in Introduction, has no solution in $\mathbb{Z}[[x]]$.
\end{remark}

\section{Fundamental solution}

Let $K$ be an arbitrary commutative unitary ring and $\frac{1}{x}K[[\frac{1}{x}]]$ the ring of formal Laurent series, those have the form
$\sum_{j=1}^{\infty}\frac{g_j}{x^j}$.
First we consider the equation
\begin{equation}\label{EQ-g}
	a_mw^{(m)}(x)+a_{m-1}w^{(m-1)}(x)+\ldots+a_1w'(x)+a_0w(x)=g(x),
\end{equation}
where $a_0,\ldots,a_m \in K$ and $g(x)\in \frac{1}{x}K[[\frac{1}{x}]]$.

\begin{theorem}
	Suppose $a_0$ is invertible and $g(x)=\sum_{j=1}^{\infty} {\frac{g_j}{x^j}}$ and the sequence $\{c_k\}$ is uniquely determined by \eqref{coeffs}. Then the series 
	\begin{equation}\label{conv-sol}
		w(x)=\sum_{k=0}^{\infty} {c_kg^{(k)}(x)},
	\end{equation}
	which is a well-defined Laurent formal series, is a unique solution of \eqref{EQ-g} in the ring $\frac{1}{x}K[[\frac{1}{x}]]$. This series can be represented as
	\begin{equation}
		w(x)=\sum_{k=1}^{\infty} {\left(\sum_{i=0}^{k-1} (-1)^{i}c_ig_{k-i}\binom{k-1}{k-i-1}i!\right)x^{-k} }.
	\end{equation}
\end{theorem}

\begin{proof}
	First let us prove the uniqueness. Let $w(x)\in \frac{1}{x}K[[\frac{1}{x}]]$ be a solution of the homogeneous equation
	$$a_mw^{(m)}(x)+a_{m-1}w^{(m-1)}(x)+\ldots+a_1w'(x)+a_0w(x)=0.$$
	Let $-k$ denote the maximal degree of a term with a non-zero coefficient in $w(x)$.  Then the maximal degree of summands $a_mw^{(m)}(x), a_{m-1}w^{(m-1)}(x), \ldots, a_1w'(x)$ is no more than $-k-1$, therefore the coefficient of $x^{-k}$ in $w(x)$ is zero.
	
	Now let us prove that \eqref{conv-sol} is well-defined. Indeed, the maximal degree of $\{g^{(k)}\}$ decreases. It follows that each degree of $x$ has finite sum of the elements of $K$ as a coefficient in the sum \eqref{conv-sol}.
	
	Similarly as for Theorem \ref{uniq-ex-abs} one can check that the series \eqref{conv-sol} is a solution of \eqref{EQ-g}.
\end{proof}

\begin{corollary}
	Let $a_0$ is invertible. Then Equation \eqref{EQ-g} with non-homogeneity $g(x)=\frac{1}{x}$ has a unique solution in the ring $\frac{1}{x}K[[\frac{1}{x}]]$
	\begin{equation}\label{F}
		\mathcal{E}(x)=\sum_{k=0}^{\infty} {c_k\frac{(-1)^kk!}{x^{k+1}}}.
	\end{equation}
\end{corollary}

\begin{example}\label{ex_2_E}
	In the case $m=2$ the sequence $\{c_k\}$ is calculated in Remark \ref{low-order-c}. Then the equation
	$$a_2w''(x)+a_1w'(x)+a_0w(x)=\frac{1}{x}$$
	has the unique solution from $\frac{1}{x}K[[\frac{1}{x}]]$:
	$$\mathcal{E}(x)=\sum_{k=0}^{\infty} {c_k\frac{(-1)^kk!}{x^{k+1}}}=\sum_{k=0}^{\infty} {\sum_{j=0}^{[\frac{k}{2}]} { {(-1)^{j}\binom{k-j}{j}a_0^{-k+j-1}a_1^{k-2j}a_2^j}}\frac{k!}{x^{k+1}}}.$$
\end{example}

In \cite{Fund-sol} we consider the convolution (Hurwitz product) of two formal Laurent series $f(x)=\sum_{i=1}^{\infty} \frac{f_i}{x^i} \in \frac{1}{x}K[[\frac{1}{x}]]$ and $g(x)=\sum_{i=1}^{\infty} \frac{g_i}{x^i} \in \frac{1}{x}K[[\frac{1}{x}]]$ be giving by 
$$(f*g)(x)=\sum_{i=0}^{\infty} {(-1)^if_{i+1}\frac{g_n^{(i)}(x)}{i!}}.$$
	
It is easy to check that $(\frac{1}{x}*g)(x)=g(x)$ for any $g\in\frac{1}{x}K[[\frac{1}{x}]]$ and $(\mathcal{E}*g)(x)$ coincides with the right side of the equality \eqref{conv-sol}, i.e. it is a unique solution of \eqref{EQ-g}, where $\mathcal{E}(x)$ is defined by \eqref{F}. For the first-order equation it is shown in \cite[Theorem 4.2]{Fund-sol}. That allows us to regard $\mathcal{E}(x)$ as a \textit{fundamental solution} of Equation \eqref{EQ-g} (see, for example, \cite[Ch.3]{Vladimirov}). 

\begin{example}
	In the case $m=2$ the formal Laurent series $\mathcal{E}(x)$ from Example \ref{ex_2_E} is a fundamental solution of the equation $a_2w''(x)+a_1w'(x)+a_0w(x)=g(x)$ in the ring $\frac{1}{x}K[[\frac{1}{x}]]$.
	
	
\end{example}

Also in \cite{Fund-sol} the convolution of a formal Laurent series $g(x)\in\frac{1}{x}K[[\frac{1}{x}]]$ with the polynomial $f(x)\in K[[x]]$ is defined as $(g*f)(x)=\sum_{i=0}^{\infty} {(-1)^ig_{i+1}\frac{f_n^{(i)}(x)}{i!}}$. It is not difficult to check that it is a unique polynomial solution of Equation \eqref{EQ}.

Now let $(K,|\cdot|)$ be a valuation ring of a field $F$, where $F$ is a complete non-Archemedean field of characteristics zero (see Section 4).  Let us consider the topology of coefficientwise convergence on the ring $K[[x]]$.

Using series \eqref{F}, one can represent the solution of Equation \eqref{EQ} as a convolution, as it is done in \cite{Hurwitz} for the first-order equation and the case $(K,|\cdot|)=(\mathbb{Z}_p,|\cdot|_p)$.

\begin{lemma}\label{well-def-conv} Suppose $b(x)=\sum_{i=1}^{\infty} \frac{b_i}{x^i}$ and $f(x)=\sum_{i=0}^{\infty} f_ix^i$. If $|b_i| \rightarrow 0$, then the series
	\begin{equation}\label{conv}
		\sum_{i=0}^{\infty} {(-1)^ib_{i+1}\frac{f^{(i)}(x)}{i!}}
	\end{equation}
converges in the topology of coefficientwise convergence on $K[[x]]$.
\end{lemma}

\begin{proof}
	Coefficients of $\frac{f^{(i)}(x)}{i!}$ belong to $K$, therefore they are no more then 1 by the valuation $|\cdot|$. Since $b_i\rightarrow 0$, then the coefficient of each degree of $x$ in \eqref{conv} is a convergent series.
\end{proof}

This lemma implies that the following notion of the convolution is well-defined.

\begin{definition}
	 Similarly as in \cite{infinite_ord} and \cite{Hurwitz}, the convolution $(b*f)(x)$ of a formal Laurent series $b(x)=\sum_{i=1}^{\infty} \frac{b_i}{x^i}\in\frac{1}{x}K[[\frac{1}{x}]]$ in that $b_i$ tends to $0$ and an arbitrary formal power series $f(x)=\sum_{i=0}^{\infty} f_ix^i\in K[[x]]$ is giving by \eqref{conv}.
\end{definition}	

\begin{lemma}[Properties of the convolution]\label{conv-prop} For any $b\in \frac{1}{x}[[\frac{1}{x}]]$ and $f\in K[[x]]$ the following equalities hold:
	\begin{enumerate}
		\item $(b*f)'(x)=(b*f')(x)=(b'*f)(x)$;
		\item $(\frac{1}{x}*f)(x)=f(x)$.
	\end{enumerate}
\end{lemma}

\begin{proof} Since $$b'(x)=-\sum_{i=1}^{\infty} \frac{ib_{i}}{x^{i+1}}\ \mbox{, then}\ (b'*f)(x)=\sum_{i=1}^{\infty} {\frac{(-1)^{i-1}b_{i}}{(i-1)!}f^{(i)}(x)}.$$
	Therefore, 	
	\begin{multline*}
		(b*f)'(x)=\sum_{i=0}^{\infty} {\frac{(-1)^ib_{i+1}}{i!}f^{(i+1)}(x)}=\\
		=(b*f')(x)= \sum_{i=1}^{\infty} {\frac{(-1)^{i-1}b_{i}}{(i-1)!}f^{(i)}(x)}=(b'*f)(x).
	\end{multline*}
\end{proof}

\begin{theorem}\label{fund-sol}
	Suppose the assumptions of Theorem \ref{exist-nonarch} hold. Then the unique solution from $K[[x]]$ of Equation \eqref{EQ} has the form
	$$w(x)=(\mathcal{E}*f)(x),$$
	where $\mathcal{E}(x)$ is defined by \eqref{F}. 
\end{theorem}

\begin{proof}
	Under the assumptions of Theorem \ref{exist-nonarch}, $\{k!c_k\}$ tends to $0$, since $\{c_k\}$ tends to $0$ (see inequality \eqref{ineq}), therefore by Lemma \ref{well-def-conv} the convolution $\mathcal{E}*f$ is well-defined.
	
	Let us denote the differential operator on the ring $\frac{1}{x}K[[\frac{1}{x}]]$: $$P(w)=a_mw^{(m)}+a_{m-1}w^{(m-1)}+\ldots+a_1w'+a_0w.$$ Then Equation \eqref{EQ-g} has the form $P(w)=f$. The solution of the equation $(P(w))(x)=\frac{1}{x}$ is given by formula \eqref{F}. 
	
	It follows that $(\mathcal{E}*f)(x)$ is a solution of $\eqref{EQ}$. Indeed, due to the properties of the convolution (Lemma \ref{conv-prop}) we get $$P(\mathcal{E}*f)(x)=(P(\mathcal{E})*f)(x)=(\frac{1}{x}*f)(x)=f(x).$$
\end{proof}

The following corollary specifies the previous result to the case of the integer coefficients.

\begin{corollary}
	Suppose $a_i \in \mathbb{Z}$, $i=0,\ldots,m$, $p$ is a prime and $p$ is not a divisor of $a_0$. Then the unique solution from $\mathbb{Z}_p[[x]]$ of Equation \eqref{EQ} has the form
	$$w(x)=(\mathcal{E}*f)(x).$$
\end{corollary}

\begin{proof}
	Suppose now $(F,|\cdot|)=(\mathbb{Q}_p,|\cdot|_p)$. Repeat the proof of Theorem \ref{fund-sol} except the proposition that $\{c_k\}$ tends to zero. Now $\{k!c_k\}$ tends to $0$ because $|k!|_p\rightarrow 0$ and $|c_k|_p\leq 1$.
\end{proof}

\begin{remark}
	Theorem \ref{fund-sol} allows us to regard the series $\mathcal{E}(x)=\sum_{k=0}^{\infty} {c_k\frac{(-1)^kk!}{x^{k+1}}}$ as a fundamental solution to Equation \eqref{EQ} for the ring $K[[x]]$.
\end{remark}

\bibliographystyle{plainurl}
\bibliography{mybibfile}

\end{document}